\documentclass[reqno]{amsart}

\usepackage{amsmath,amsthm,amssymb,amsfonts,mathrsfs,amscd}
\usepackage{multirow}
\usepackage{graphicx}
\usepackage{bm}
\usepackage{mathrsfs}
\usepackage{dsfont}
\usepackage{tikz}
\usepackage{tikz-3dplot}
\usepackage{subfig}
\usepackage{stmaryrd}

\newtheorem{remark}{Remark}[section]

\newtheorem{lemma}{Lemma}[section]
\newtheorem{theorem}{Theorem}[section]

\makeatletter
  \@addtoreset{equation}{section}
  \newcommand\figcaption{\def\@captype{figure}\caption}
  \newcommand\tabcaption{\def\@captype{table}\caption}
\makeatletter

\begin{document}

\title{A variant of the plane wave least squares method for the time-harmonic Maxwell's equations}

\author{Qiya Hu}
\author{Rongrong Song}

\thanks{1. LSEC, ICMSEC, Academy of Mathematics and Systems Science, Chinese Academy of Sciences, Beijing
100190, China; 2. School of Mathematical Sciences, University of Chinese Academy of Sciences, Beijing 100049,
China (hqy@lsec.cc.ac.cn, songrongrong@lsec.cc.ac.cn). This work was funded by Natural Science Foundation of China G11571352.}

\maketitle

{\bf Abstract.}
In this paper we are concerned with the plane wave method for the discretization of time-harmonic Maxwell's equations in three dimensions. As pointed out in \cite{pwdg}, it is difficult to derive a satisfactory $L^2$ error estimate of the standard plane wave approximation of the time-harmonic Maxwell's equations. We propose a variant of the plane wave least squares (PWLS) method and show that the new plane wave approximations possess the desired $L^2$ error estimate.
Moreover, the numerical results indicate that the new approximations have sightly smaller $L^2$ errors than the standard plane wave approximations. More importantly, the results are derived for more general models in layered media.



\pagestyle{myheadings}
\thispagestyle{plain}
\markboth{}{QIYA HU AND RONGRONG SONG}

\section{Introduction}
The plane wave method was first introduced to solve Helmholtz equations and was then extended to solve time-harmonic Maxwell's equations. This method is different from the traditional finite element method because of its special choice of basis functions. In the plane wave methods, the basis functions are chosen as the exact solutions of the differential equations without boundary conditions, so the resulting approximate solutions possess higher accuracies than that generated by the other methods. Over the last ten years, various plane wave methods were proposed in the literature, for example, the ultra weak variational formulation (UWVF) \cite{Buffa2008,ref11,hmm}, the plane wave discontinuous Galerkin (PWDG) method \cite{Hiptmair2011,pwdg} and the plane wave least squares (PWLS) method \cite{Hu2014,long,Hu2017,ref12} (more references on the plane wave methods can be found in the survey article \cite{Hiptmair2015}).

Although there are many works to introduce plane wave methods and analyze the convergence of plane wave approximations for Helmholtz equation and time-harmonic Maxwell's equations, a satisfactory $L^2$ error estimate of the plane wave
approximations for time-harmonic Maxwell's equations has not been obtained yet (see \cite{pwdg} for more detailed explanations on the reasons).
Moreover, to our knowledge, there is no convergence analysis on the plane wave methods for Helmholtz equation or time-harmonic Maxwell's equations in layer media.



In this paper, we consider the PWLS method since it results in a Hermitian and positive definite system. We propose a variant of the PWLS method to discretize time-harmonic homogeneous Maxwell's
equations in three dimensions. In the new PWLS method, we add the jump of the normal complement of the solutions into the objective functional defined in the existing PWLS methods. We prove that
the errors of the resulting plane wave approximations possess a desired $L^2$ estimate that is almost the same as the one in the case of Helmholtz equation.
Here we need not to assume that the material coefficients are constants, namely, layer media are allowed.
The numerical results show that the plane wave approximations generated by the proposed method is slightly more accurate than that generated by the existing PWLS methods.

The paper is organized as follows:
In Section 2, we give a description of the second order system of Maxwell equations associated with
triangulations.
In Section 3, we present the proposed PWLS method for time-harmonic Maxwell's equations.
In Section 4, we explain how to discretize the variational problem.
In Section 5, we give a $L^2$ error estimate for the approximate solutions generated by the new variational problem.
In Section 6, we simplify the variational formulation for the case of homogeneous material.
Finally, in Section 7 we report some numerical results to confirm the effectiveness of the new method.

\section{Description of time-harmonic Maxwell's equations}

In this section we recall the first order system of Maxwell's equations and derive the corresponding second order system.

Let $\Omega \subset \mathbb{R}^3$ be a bounded, polyhedral domain. We denote by ${\bf n}$ the unit normal vector field on $\partial\Omega$ pointing outside $\Omega$.
We consider the following formulation of the three-dimensional time-harmonic Maxwell's equations in terms of electric filed ${\bf E}$ and magnetic field ${\bf H}$ written as the first-order system of equations
in $\Omega$
\begin{equation} \label{eq2 1}
\left\{
\begin{aligned}
& \nabla\times{\bf E} - i\omega\mu{\bf H}=0, \\
& \nabla\times{\bf H} + i\omega\varepsilon{\bf E}=0, \\
& \nabla\cdot(\varepsilon{\bf E})=0,\\
& \nabla\cdot(\mu{\bf H})=0,
\end{aligned}
\right.
\end{equation}
with the lowest-order absorbing boundary condition
\begin{equation} \label{eq2 2}
-{\bf E}\times{\bf n}+\sigma({\bf H}\times{\bf n})\times{\bf n}={\bf g}
\quad \text{on} \ \gamma=\partial\Omega.
\end{equation}

Here $\omega>0$ denotes the temporal frequency of the field, and ${\bf g}\in L^2(\partial\Omega)$. The material coefficients $\varepsilon,\mu$ and $\sigma$ are understood as usual (refer to \cite{hmm}). Notice that we often have $\sigma=\sqrt{\mu/|\varepsilon|}$ in applications.

In particular, if $\varepsilon$ is complex valued, then the material is known as an absorbing medium; otherwise the material is called a non-absorbing medium. And $\varepsilon$ may have different values in different areas, but if $\varepsilon$ is a constant in the whole domain, then the medium is called homogeneous medium.

Based on the first equation of (\ref{eq2 1}), we can write ${\bf H}$ in terms of ${\bf E}$ as ${\bf H}=\frac{1}{i\omega\mu} \nabla\times{\bf E}$. Substituting this expression into the second equation and into the boundary condition, we obtain the second-order Maxwell's equations
\begin{equation} \label{eq2 3}
\left\{
\begin{aligned}
& \nabla\times(\frac{1}{i\omega\mu}\nabla\times{\bf E}) + i\omega\varepsilon {\bf E}=0
& \text{in}\ \Omega,\\
& -{\bf E}\times{\bf n}+\frac{\sigma}{i\omega\mu}((\nabla\times{\bf E})\times{\bf n})\times{\bf n}={\bf g}
& \text{on}\ \gamma.
\end{aligned}
\right.
\end{equation}

\section{Variational formulation for Maxwell's equations}

In this section we introduce details of the new variational problem of the Maxwell's equations based on triangulation.

Let $\Omega$ be divided into a partition in the sense that
$$ \overline{\Omega}=\bigcup_{k=1}^N\overline{\Omega}_k,\quad \Omega_l\bigcap\Omega_j=\emptyset, \quad\text{for } l\not=j, $$
and let ${\mathcal T}_h$ denote the triangulation comprised of elements $\{\Omega_k\}_{k=1}^N$, where $h$ denotes the meshwidth of the triangulation.
Define $$ \Gamma_{lj}=\partial\Omega_l\bigcap\partial\Omega_j, \quad\text{for } l\not=j, $$
and $$ \gamma_k=\partial\Omega_k \bigcap \partial\Omega, \quad k=1,\ldots,N, \quad\gamma=\bigcup^N_{k=1}\gamma_k. $$

For each element $\Omega_k$, let ${\bf E}|_{\Omega_k}={\bf E}_k$ and $\varepsilon|_{\Omega_k}=\varepsilon_k$. It is known that $\omega$ is a constant on $\Omega$.
As usual, we can assume that $\mu$ and $\varepsilon$ are constants on each element. If $\varepsilon_j\not=\varepsilon_k$ for two different $k$ and $j$, the media on $\Omega$
is called layer media.

Notice that we have the hidden relation $div(\varepsilon {\bf E})=0$ from the first equation of (\ref{eq2 3}).
Then the reference problem (\ref{eq2 3}) will be solved by finding the local electric field ${\bf E}_k$ such that
\begin{equation}
\nabla\times(\frac{1}{i\omega\mu}\nabla\times{\bf E}_k) + i\omega\varepsilon_k {\bf E}_k=0 \quad \text{in} \ \Omega_k,    \label{eq2 4}
\end{equation}
with the transmission conditions on each interface $\Gamma_{lj}$ (notice that ${\bf n}_l=-{\bf n}_j \ \text{on}~\Gamma_{lj}$)
\begin{eqnarray}
\left\{
\begin{array}{ll}
{\bf E}_l\times{\bf n}_l+{\bf E}_j\times{\bf n}_j=0,\\
(\frac{1}{i\omega\mu}\nabla\times{\bf E}_l)\times{\bf n}_l+(\frac{1}{i\omega\mu}\nabla\times{\bf E}_j)\times{\bf n}_j=0, \\
\varepsilon_l{\bf E}_l\cdot{\bf n}_l+\varepsilon_j{\bf E}_j\cdot{\bf n}_j=0. \\
\end{array}
\right.
\label{eq2 5}
\end{eqnarray}

The boundary condition becomes
\begin{equation}\label{eq2 6}
-{\bf E}_k\times{\bf n}+\frac{\sigma}{i\omega\mu}((\nabla\times{\bf E}_k)\times{\bf n})\times{\bf n}={\bf g}, \quad \text{on} \ \gamma_k.
\end{equation}

For an element $\Omega_k$, let $\mbox{{\bf H}}(\text{curl},\Omega_k)$ denote the standard Sobolev space. Set
\begin{equation} \label{eq3 1}
{\bf V}(\Omega_k)=\{{\bf E}_k\in \mbox{{\bf H}}(\text{curl},\Omega_k), ~{\bf E}_k \text{~satisfies the equation (\ref{eq2 4})}\}\,
\end{equation}
and define
$$ {\bf V}({\mathcal T}_h)=\prod\limits_{k=1}^N {\bf V}(\Omega_k) $$

For ${\bf F}\in {\bf V}({\mathcal T}_h)$, set ${\bf F}|_{\Omega_k}={\bf F}_k$. For ease of notation, define
$$ \Phi({\bf F_k})=\frac{\sigma}{i\omega\mu}((\nabla\times{\bf F_k})\times{\bf n_k})  $$ 
and
$$ \Psi({\bf F}_k)={1\over i\omega\mu}(\nabla\times{\bf F}_k). $$ 

Then the boundary condition(\ref{eq2 6}) can be written as
\begin{equation}
\label{eq3 2} -{\bf F}_k\times{\bf n}_k+\Phi({\bf F}_k)\times{\bf n}_k={\bf g} \quad \text{on} ~\gamma_k.
\end{equation}

For each local interface $\Gamma_{lj}$ $(l<j)$, we define the jumps on $\Gamma_{lj}$ as follows (note that ${\bf n}_l=-{\bf n}_j$)
\begin{eqnarray} \label{eq3 3}
\left\{
\begin{array}{ll}
\llbracket{\bf F}\times{\bf n}\rrbracket = {\bf F}_l\times{\bf n}_l+{\bf F}_j\times{\bf n}_j, \\
\llbracket\Psi({\bf F})\times{\bf n}\rrbracket = \Psi({\bf F}_l)\times{\bf n}_l + \Psi({\bf F}_j)\times{\bf n}_j,  \\
\llbracket\varepsilon{\bf F}\cdot{\bf n}\rrbracket = \varepsilon_l{\bf F}_l\times{\bf n}_l +\varepsilon_j{\bf F}_j\cdot{\bf n}_j. \\
\end{array}
\right.
\end{eqnarray}

It is easy to see that the transmission conditions (\ref{eq2 5}) are equivalent to the following conditions
\begin{equation} \label{eq3 4}
\llbracket{\bf F}\times{\bf n}\rrbracket=0, ~~~
\llbracket\Psi({\bf F})\times{\bf n}\rrbracket=0 ~~~\mbox{and}~~~
\llbracket\varepsilon{\bf F}\cdot{\bf n}\rrbracket=0 \quad
\text{on} ~\Gamma_{lj}.
\end{equation}

A vector field ${\bf E}\in {\bf V}({\mathcal T}_h)$ is the solution of (\ref{eq2 4}), (\ref{eq2 5}) and (\ref{eq2 6}) if and only if the conditions (\ref{eq3 2}) and (\ref{eq3 4}) are satisfied.
Based on this observation, we define the functional
\begin{equation}  \label{minfun}
\begin{split}
J({\bf F})=
& \sum\limits_{k=1}^N \delta\int_{\gamma_k}|-{\bf F}_k\times{\bf n}_k + \Phi({\bf F}_k)\times{\bf n}_k-{\bf g}|^2 ~ds \\
& +\sum\limits_{l<j} \Big(\alpha\int_{\Gamma_{lj}}|\llbracket{\bf F}\times{\bf n} \rrbracket|^2 ~ds +\beta\int_{\Gamma_{lj}}|\llbracket\Psi({\bf F})\times{\bf n} \rrbracket|^2 ~ds \\
& +\theta\int_{\Gamma_{lj}}|\llbracket\varepsilon{\bf F}\cdot{\bf n} \rrbracket|^2 ~ds \Big)
\quad \forall{\bf F}\in {\bf V}({\mathcal T}_h),
\end{split}
\end{equation}
with $\delta,~\alpha,~\beta$ and $\theta$ being positive numbers, which will be specified later. It is clear that $J({\bf F})\geq 0$. Then we consider the following minimization problem:
Find ${\bf E}\in {\bf V}({\mathcal T}_h)$ such that
\begin{equation}
J({\bf E})=\min\limits_{{\bf F}\in {\bf V}({\mathcal T}_h)}J({\bf F}).\label{min}
\end{equation}

If ${\bf E}$ is the solution of the equations (\ref{eq2 4}), (\ref{eq2 5}) and (\ref{eq2 6}), i.e. ${\bf E}\in {\bf V}({\mathcal T}_h)$ satisfies the conditions (\ref{eq3 2}) and (\ref{eq3 4}), then we have $J({\bf E})=0$, which implies that ${\bf E}$ is also the solution of the minimization problem (\ref{min}). When $J({\bf E})$ is very small, the vector field ${\bf E}$ should be an approximate solution of (\ref{eq2 4}), (\ref{eq2 5}) and (\ref{eq2 6}).

The variational problem associated with the minimization problem(\ref{min}) can be expressed as follows: Find ${\bf E}\in {\bf V}({\mathcal T}_h)$ such that
\begin{equation}
\begin{split}
& \sum_{k=1}^N \delta\int_{\gamma_k}(-{\bf E}_k\times{\bf n}_k +\Phi({\bf E}_k)\times{\bf n}_k -{\bf g}) \cdot\overline{-{\bf F}_k\times{\bf n}_k +\Phi({\bf F}_k)\times{\bf n}_k} ~ds \\
& +\sum_{l<j} \Big(\alpha\int_{\Gamma_{lj}}\llbracket{\bf E}\times{\bf n}\rrbracket \cdot\overline{\llbracket{\bf F}\times{\bf n} \rrbracket} ~ds +\beta\int_{\Gamma_{lj}} \llbracket\Psi({\bf E})\times{\bf n}\rrbracket\cdot\overline{\llbracket\Psi({\bf F})\times{\bf n}\rrbracket} ~ds \\
& +\theta\int_{\Gamma_{lj}}\llbracket\varepsilon{\bf E}\cdot{\bf n} \rrbracket\cdot\overline{\llbracket\varepsilon{\bf F}\cdot{\bf n} \rrbracket} ~ds \Big) =0
\quad\quad \forall ~{\bf F}\in {\bf V}({\mathcal T}_h).
\end{split}
\label{eq3 6}
\end{equation}

The method described above can be viewed as a variant of the plane wave least-squares (PWLS) method proposed in \cite{ref12} or \cite{long}, here we added the final term containing the jump of the normal complements of
the considered vector-valued functions.

Define the sesquilinear form $A(\cdot,\cdot)$ by
\begin{equation}
\begin{split}
A({\bf E},{\bf F})
& =\sum_{k=1}^N \delta\int_{\gamma_k}(-{\bf E}_k\times{\bf n}_k+\Phi({\bf E}_k)\times{\bf n}_k)\cdot\overline{-{\bf F}_k\times{\bf n}_k+\Phi({\bf F}_k)\times{\bf n}_k} ~ds \\
& +\sum_{l<j}\Big(\alpha\int_{\Gamma_{lj}}\llbracket{\bf E}\times{\bf n}\rrbracket\cdot\overline{\llbracket{\bf F}\times{\bf n} \rrbracket} ~ds +\beta\int_{\Gamma_{lj}}\llbracket\Psi({\bf E})\times{\bf n} \rrbracket\cdot\overline{ \llbracket\Psi({\bf F})\times{\bf n} \rrbracket} ~ds \\
& +\theta\int_{\Gamma_{lj}}\llbracket\varepsilon{\bf E}\cdot{\bf n}\rrbracket\cdot\overline{\llbracket\varepsilon{\bf F}\cdot{\bf n} \rrbracket} ~ds \Big)
\quad\quad \forall ~{\bf F}\in {\bf V}({\mathcal T}_h),
\end{split}
\label{eq3 7}
\end{equation}
and ${\bm \xi}\in {\bf V}({\mathcal T}_h)$, via the Riesz representation theorem, by
\begin{equation}
(\bm\xi,~{\bf F}) =
\sum_{k=1}^N \delta\int_{\gamma_k}{\bf g}\cdot\overline{-{\bf F}_k\times{\bf n}_k+\Phi({\bf F}_k)\times{\bf n}_k} ~ds  \quad\quad \forall ~{\bf F}\in {\bf V}({\mathcal T}_h).
\end{equation}

Then (\ref{eq3 6}) can be written as
\begin{eqnarray}
\left
\{\begin{array}{ll}
\text{Find} \ {\bf E}\in {\bf V}({\mathcal T}_h), \ \text{s.t.} \\
A({\bf E},{\bf F})=({\bm \xi},~{\bf F}), \quad \forall ~{\bf F}\in {\bf V}({\mathcal T}_h).
\end{array}
\right.
\label{eq3 9}
\end{eqnarray}

\begin{theorem}
Let ${\bf E}\in {\bf V}({\mathcal T}_h)$. For $k=1,\ldots,N$, assume that ${\bf E}_k\in H^{1+\delta_k}(\Omega_k)$ with $\delta_k>{1\over2}$ such that $(\nabla\times{\bf E}_k)\times{\bf n}_k\in L^2(\partial\Omega_k)$. Then the reference problem (\ref{eq2 3}) is equivalent to the new variational problem (\ref{eq3 9}).
\end{theorem}

\paragraph*{Proof.}

We only need to verify the uniqueness of the solution of problem(\ref{eq3 9}). The process of verification is standard.

Let us consider the two solutions ${\bf E}=({\bf E}_1,\ldots,{\bf E}_N)$ and ${\bf E}'=({\bf E}_1',\ldots,{\bf E}_N')$ of the variational problem, and let $\tilde{{\bf E}}=(\tilde{{\bf E}}_1,\ldots,\tilde{{\bf E}}_N)$ denote the difference between the two solutions. It follows by (\ref{eq3 9}) that the difference satisfies
\begin{equation}
A(\tilde{\bf E}, {\bf F})=0, \quad \forall ~{\bf F}\in{\bf V}({\mathcal T}_h). \label{eq3 10}
\end{equation}

Taking ${\bf F}=\tilde{{\bf E}}$, the equation (\ref{eq3 10}) becomes
\begin{equation*}
A(\tilde{\bf E},\tilde{\bf E})=0,
\end{equation*}
which implies that
$$ \int _{\gamma_k}|-\tilde{\bf E}_k\times{\bf n}_k+\Phi(\tilde{\bf E}_k)\times{\bf n}_k|^2 ~ds=0 $$
and
$$ \int_{\Gamma_{kj}}|\llbracket\tilde{\bf E}\times{\bf n} \rrbracket|^2 ~ds=\int_{\Gamma_{kj}}|\llbracket\Psi(\tilde{\bf E})\times{\bf n} \rrbracket|^2 ~ds =\int_{\Gamma_{kj}}|\llbracket\varepsilon\tilde{\bf E}\cdot{\bf n} \rrbracket|^2 ~ds = 0. $$

These show that $\tilde{{\bf E}}$ satisfies the condition (\ref{eq3 4}), i.e. the transmission continuity (\ref{eq2 5}), and verify the initial Maxwell reference problem (\ref{eq2 3}) (notice that $\tilde{{\bf E}}\in {\bf V}({\mathcal T}_h)$) with the boundary condition. Thus the function $\tilde{{\bf E}}$ vanishes on $\Omega$, which proves the uniqueness of the solution of (\ref{eq3 9}).
$\hfill \Box$

\section{Discretization of the variational problem}

In this section, we will describe a discretization of the variational problem (\ref{eq3 9}). It is based on a finite dimensional space ${\bf V}_p({\mathcal T}_h)\subset {\bf V}({\mathcal T}_h)$. We first give the precise definition of the space ${\bf V}_p({\mathcal T}_h)$.

\subsection{Basis functions of ${\bf V}_p({\mathcal T}_h)$}

In this subsection we construct basis functions with which to discretize the PWLS method.

In practice, following \cite{refOC}, a suitable family of plane waves, which are solutions of the constant-coefficient Maxwell equations, are generated on $\Omega_k$ by choosing $p$ unit propagation directions ${\bf d}_{l}(l=1,\ldots,p)$,
and defining a real unit polarization vector ${\bf G}_{l}$ orthogonal to ${\bf d}_{l}$. Then the propagation directions and polarization vectors define the complex polarization vectors ${\bf F}_{l}$ and ${\bf F}_{l+p}$ by
$$ {\bf F}_{l}={\bf G}_{l} + \text{i}{\bf G}_{l}\times {\bf d}_{l},\quad {\bf F}_{l+p}={\bf G}_{l} - \text{i}{\bf G}_{l}\times {\bf d}_{l} \quad (l=1,\ldots,p). $$

It is clear that (\ref{eq2 4}) can be written as
\begin{equation} \label{eq2 44}
\nabla\times(\nabla\times{\bf E}_k)- \kappa^2{\bf E}_k=0  \quad\text{in} ~\Omega_k,
\end{equation}
with $\kappa=\omega\sqrt{\mu\varepsilon}$. We then define the complex functions ${\bf E}_{l}$
\begin{equation}
{\bf E}_{l}=\sqrt{\mu}~{\bf F}_{l}~\text{exp}(\text{i}\kappa{\bf d}_{l}\cdot {\bf x}) \quad\quad l=1,\ldots,2p.  \label{eq4 1}
\end{equation}

It is easy to verify that every function ${\bf E}_{l}$ ($l=1,\ldots,2p$) satisfies Maxwell's equation (\ref{eq2 44}). Notice that
$\kappa$ has different values on different elements for layer media, so the plane wave basis functions ${\bf E}_{l}$ may have different
form on different elements.

Let ${\mathcal Q}_{2p}$ denote the space spanned by the $2p$ plane
wave functions ${\bf E}_{l}$ ($l=1,\ldots,2p$). Define the finite element space
\begin{equation}
{\bf V}_p({\mathcal T}_h)=\bigg\{{\bf v}\in L^2(\Omega):~{\bf v}|_K\in{\mathcal Q}_{2p} \quad \forall~K\in {\mathcal T}_h \bigg\},
\label{eq4 2}
\end{equation}
which has $N\times 2p$ basis functions, which are defined by
\begin{eqnarray}
\bm\phi^k_{l}({\bf x})=
\left\{
\begin{array}{ll}
{\bf E}_{l}({\bf x}) \quad {\bf x}\in\Omega_k \\
0 \quad\quad\quad {\bf x} \in \Omega_j \mbox{ when } j\neq k
\end{array}
\right.
~~(k=1,\ldots,N, ~l=1,\ldots,2p).
\label{eq4 3}
\end{eqnarray}

\subsection{Discrete problem and algebraic form of (\ref{eq3 9})}

Let ${\bf V}_p({\mathcal T}_h)$ be defined as in the previous subsection. Now we define an approximation of ${\bf E}_k$ by
\begin{equation}
{\bf E}^k_h=\sum_{l=1}^{2p}x_l^k \bm{\phi}^k_{l}.
\end{equation}

In this case, the discrete variational problem associated with (\ref{eq3 9}) can be described as follows
\begin{eqnarray}
\left\{
\begin{array}{ll}
\text{Find} \ {\bf E}_h\in {\bf V}_p({\mathcal T}_h), \ \text{ s.t.} \\
A({\bf E}_h,{\bf F}_h)=(\bm\xi,~{\bf F}_h),~~\forall ~{\bf F}_h\in {\bf V}_p({\mathcal T}_h).
\end{array}
\right.
\label{eq4 5}
\end{eqnarray}

After solving problem (\ref{eq4 5}), the approximate solution of Maxwell's equations (\ref{eq2 3}) are obtained directly, because the unknown ${\bf E}_h$ is defined on the elements. Moreover, the structure of the sesquilinear form $A(\cdot, \cdot)$ is very simple, so the method seems easier to implement.

As in \cite{long}, let ${\mathcal A}$ be the stiffness matrix associated with the sesquilinear form $A(\cdot,\cdot)$, and let $b$ denote the vector associated with the vector product $(\bm\xi,~{\bf F}_h)$. Then the discretized problem (\ref{eq4 5}) leads to the algebraic system
\begin{equation}
{\mathcal A}X=b, \label{eq4 6}
\end{equation}
where $ X=(x_{11},x_{12},\ldots,x_{1\,2p},x_{21},\ldots,x_{2\,2p},\ldots,x_{N1},\ldots,x_{N\,2p})^t\in \mathbb{C}^{2pN} $ is the unknown vector. Furthermore, we know that the matrix ${\mathcal A}$ is Hermitian positive definite from the definition of the bilinear form $A(\cdot,\cdot)$, so the system (\ref{eq4 6}) is easier to be solved.

\section{$L^2$ Error estimate of the plane wave approximations}

In this section, we derive the $L^2$ estimate for the errors of the approximate solution ${\bf E}_h$ defined by (\ref{eq4 5}). The analysis is based on the techniques developed in \cite{pwdg}, as well as some new observations.

Set ${\bf e}_h={\bf E}-{\bf E}_h$, where ${\bf E}$ and ${\bf E}_h$ are defined by (\ref{eq3 9}) and (\ref{eq4 5}) respectively. Clearly we have
\begin{equation} \label{eq1}
\nabla\times(\nabla\times{{\bf e}_h})- \kappa^2{{\bf e}_h}=0   \quad \quad\text{in}~ \Omega_k.
\end{equation}

It is clear that $A({\bf F},{\bf F}) \ge 0$ and that $A({\bf F},{\bf F})=0$ for ${\bf F}\in {\bf V}({\mathcal T}_h)$ if and only if ${\bf F}=0$. Thus, $A(\cdot,\cdot)$ is a norm on ${\bf V}({\mathcal T}_h)$. For ease of notation, this norm is denoted by $||\cdot||_{\bf V}$ in the rest of this section.

The key technique is to use a Helmholtz-type decomposition of ${\bf e}_h$
\begin{equation} \label{eq2}
{\bf e}_h = {\bf w} + \nabla p,
\end{equation}
where $p\in H_0^1(\Omega)$ is defined by
\begin{equation}
\int_\Omega (\overline{\varepsilon}\nabla p)\cdot \nabla\overline{q} ~d\text{x} = \int_\Omega (\overline{\varepsilon}{\bf e}_h) \cdot \nabla\overline{q} ~d\text{x},\quad
\forall \, q \in H_0^1(\Omega)\label{5.new1}
\end{equation}
and ${\bf w}\in H(\text{curl},\Omega)$ satisfying $div(\overline{\varepsilon}{\bf w})=0$. We point out that, for the case of layer media (i.e., $\varepsilon$ is not a constant),
the definitions of the functions $p$ and ${\bf w}$ in the decomposition (\ref{eq2}) are different from that in the Helmholtz decomposition considered in \cite{pwdg}. 

Throughout this
paper, we always use $|\varepsilon|_\text{max}$ (rep. $|\varepsilon|_\text{min}$) to denote the maximal (rep. minimal) value of the function $|\varepsilon(x)|$ on $\bar{\Omega}$.  

\begin{lemma} \label{wpes}
For ${\bf w}$ and $p$ defined above, the stabilities hold
\begin{equation}
||{\bf w}||_{0,\Omega} \leq \frac{|\varepsilon|_\text{max}}{|\varepsilon|_\text{min}}||{\bf e}_h||_{0,\Omega}
\quad\mbox{and}\quad
||p||_{1,\Omega} \leq |p|_{1,\Omega} \leq \frac{|\varepsilon|_\text{max}}{|\varepsilon|_\text{min}}||{\bf e}_h||_{0,\Omega}.
\end{equation}
\end{lemma}

\paragraph*{Proof}
From the equality $div(\overline{\varepsilon}{\bf w})=0$, we have
$$\int_\Omega \overline{\varepsilon}{\bf w}\cdot \overline{\nabla p} ~d\text{x}
=-\int_\Omega div(\overline{\varepsilon}{\bf w})\cdot \overline{p} ~d\text{x} +\int_{\partial \Omega} (\overline{\varepsilon}{\bf w}\cdot{\bf n})\cdot \overline{p} ~ds = 0. $$

Then, by (\ref{eq2}) and (\ref{5.new1}) we can easily get
\begin{equation*}
\begin{split}
|\varepsilon|_\text{min} ||{\bf w}||_{0,\Omega}^2
& \leq |\int_\Omega \overline{\varepsilon}{\bf w} \cdot \overline{\bf w} ~d\text{x}|
=|\int_\Omega \overline{\varepsilon}{\bf w} \cdot \overline{({\bf e}_h-\nabla p)} ~d\text{x}| \\
& =|\int_\Omega \overline{\varepsilon}{\bf w} \cdot \overline{\bf e}_h ~d\text{x}|
\leq |\varepsilon|_\text{max}||{\bf w}||_{0,\Omega} ||{\bf e}_h||_{0,\Omega},
\end{split}
\end{equation*}
and
\begin{equation*}
\begin{split}
|\varepsilon|_\text{min} ||\nabla p||_{0,\Omega}^2
& \leq |\int_\Omega \varepsilon \nabla p \cdot \overline{\nabla p} ~d\text{x}|
=|\int_\Omega \varepsilon \nabla p \cdot \overline{({\bf e}_h-{\bf w})} ~d\text{x}| \\
& =|\int_\Omega \varepsilon \nabla p \cdot \overline{\bf e}_h ~d\text{x}|
\leq |\varepsilon|_\text{max}||\nabla p||_{0,\Omega} ||{\bf e}_h||_{0,\Omega}.
\end{split}
\end{equation*}

Thus
\begin{equation*}
||{\bf w}||_{0,\Omega}
\leq \frac{|\varepsilon|_\text{max}}{|\varepsilon|_\text{min}} ||{\bf e}_h||_{0,\Omega}
\quad\mbox{and}\quad
||p||_{1,\Omega}
\leq |p|_{1,\Omega}
\leq \frac{|\varepsilon|_\text{max}}{|\varepsilon|_\text{min}} ||{\bf e}_h||_{0,\Omega}.
\end{equation*}
$\hfill \Box$

In order to develop the argument needed for the error analysis below, we will consider the adjoint Maxwell's equations (the first equation is different from that in \cite{pwdg})
\begin{equation} \label{adj}
\left\{
\begin{aligned}
& \nabla\times(\nabla\times {\bf u})- \overline{\kappa}^2 {\bf u} = \overline{\varepsilon}{\bf w}  \quad\quad \text{in}\ \Omega,\\
& {\bf u}\times {\bf n} + \frac{\sigma}{i\omega\mu}\big((\nabla \times {\bf u})\times{\bf n}\big)\times{\bf n} = 0 \quad\quad \text{on} \ \partial\Omega.
\end{aligned}
\right.
\end{equation}

In the rest of this paper, we always use $C$ to denote a generic positive constant independent of $h$, $p$ and $\omega$, but its value might change at different occurrence. For ease of notation,
we choose $\sigma=\sqrt{\mu/|\varepsilon|}$.

\begin{lemma}\cite{stab} \label{stability}Under the previous assumptions on $\Omega$, the solution $\bf u$ of problem (\ref{adj}) belongs to $H^{1/2+s}(\text{curl},\,\Omega)$, for all the real parameters $s>0$ such that $s\leq \tilde s$, where $0 < \tilde s < 1/2$ is a parameter only depending on $\Omega$. Moreover, there is positive constant $C$ independent of $\omega$ and $\varepsilon$, such that
\begin{equation} \label{esti}
||\nabla \times {\bf u}||_{0,\Omega} +\omega ||\overline{\varepsilon}^{1/2}{\bf u}||_{0,\Omega} \leq C||\overline{\varepsilon}{\bf w}||_{0,\Omega}
\end{equation}
and
\begin{equation}  \label{estiesti}
||\nabla \times {\bf u}||_{1/2+s,\Omega} +\omega |\varepsilon|^{1/2}_\text{min}||{\bf u}||_{1/2+s,\Omega} \leq C\,\omega|\varepsilon|^{1/2}_\text{max} ||\overline{\varepsilon}{\bf w}||_{0,\Omega}.
\end{equation}
When $\Omega$ is convex, the inequality (\ref{estiesti}) holds for all $0 < s < 1/2$ .
\end{lemma}
$\hfill \Box$

\begin{remark} We would like to give simple explanations to above estimates. From the analysis in \cite[Thm.3.1]{stab}, we can get
\begin{equation*}
||\nabla \times {\bf u}||_{0,\Omega} +\omega ||\overline{\varepsilon}^{1/2}{\bf u}||_{0,\Omega}
\!\leq\! \big( \frac{d^2}{\gamma} (|\varepsilon|^\frac{1}{2}_\text{max}|\sigma|_\text{min} +4\mu|\varepsilon|^{-\frac{1}{2}}_\text{max}|\sigma|^{-1}_\text{min}) +d^\frac{1}{2}\mu^\frac{1}{2} \big)||\overline{\varepsilon}{\bf w}||_{0,\Omega}.
\end{equation*}
Then (\ref{esti}) can be derived directly by this inequality since $\sigma=\sqrt{\mu/|\varepsilon|}$. Furthermore, from the analysis in \cite[Thm.4.4]{stab}, we can obtain (\ref{estiesti}). It is not hard to see that the constant $C$ in Lemma \ref{stability} is independent of the parameters $\omega$ and $\varepsilon$.
\end{remark}


\begin{lemma} \label{LL} Let the parameters $\delta,~\alpha,~\beta$ and $\theta$ in the functional (\ref{minfun}) be chosen as
\begin{equation} \label{parameters}
\begin{split}
\delta=\alpha=\frac{|\varepsilon|^4_\text{max}}{|\varepsilon|^4_\text{min}},\,\quad \beta=\frac{|\varepsilon|^4_\text{max}}{|\varepsilon|^5_\text{min}} \quad\mbox{ and}\quad \theta=\frac{|\varepsilon|^2_\text{max}}{|\varepsilon|^4_\text{min}}.
\end{split}
\end{equation}
Then the error ${\bf e}_h={\bf E}-{\bf E}_h$ satisfies
\begin{equation} \label{most}
||{\bf e}_h||_{0,\Omega}
\leq C \omega^{1/2}h^{-1/2}||{\bf e}_h||_{\bf V}.
\end{equation}
\end{lemma}

\paragraph*{Proof.}Using the decomposition (\ref{eq2}), yields
\begin{equation}
|\varepsilon|_\text{min}||{\bf e}_h||^2_{0,\Omega}
\leq |\int_\Omega \varepsilon{\bf e}_h \cdot \overline{\bf e}_h ~d\text{x}|
=|\int_\Omega {\bf e}_h \cdot \overline{\overline{\varepsilon}{\bf w}} ~d\text{x}+\int_\Omega \varepsilon{\bf e}_h \cdot \overline{ \nabla p} ~d\text{x}|.\label{5.new4}
\end{equation}

Integrating by parts, using equality (\ref{eq1}) and taking into account the boundary conditions of (\ref{adj}), we deduce that
\begin{equation}
\begin{split}
\int_\Omega\! {\bf e}_h \!\cdot\! \overline{\overline{\varepsilon}{\bf w}} ~d\text{x}
& = \sum_{k=1}^N \int_{\Omega_k}{\bf e}_h \cdot \overline{\big(\nabla\times(\nabla\times {\bf u})- \overline{\kappa}^2 {\bf u}\big)} ~d\text{x} \\
& = \sum_{k=1}^N \bigg(\int_{\Omega_k} \big(\nabla\times(\nabla\times {\bf e}_h)\!-\!\kappa^2{\bf e}_h\big)\cdot\overline{\bf u} ~d\text{x} \!-\!\int_{\partial \Omega_k}{\bf e}_h\cdot\overline{(\nabla\times{\bf u})\times{\bf n_k}} ~ds \\
& \quad -\int_{\partial \Omega_k}(\nabla\times{\bf e}_h)\cdot\overline{({\bf u}\times{\bf n_k})} ~ds \bigg)\\
& = \sum_{k=1}^N \bigg( \int_{\partial \Omega_k}({\bf e}_h\times{\bf n_k})\cdot\overline{(\nabla\times{\bf u})} ~ds + \int_{\partial \Omega_k}(\nabla\times{\bf e}_h)\times{\bf n_k}\cdot\overline{\bf u} ~ds \bigg)\\
& = \sum_{l<j} \bigg( \int_{\Gamma_{lj}}\llbracket{\bf e}_h\times{\bf n}\rrbracket\cdot\overline{(\nabla\times{\bf u})} ~ds + \int_{\Gamma_{lj}}\llbracket(\nabla\times{\bf e}_h)\times{\bf n}\rrbracket\cdot\overline{\bf u} ~ds \bigg) \\
& \quad +\sum_{k=1}^N \int_{\gamma_k}(-{\bf e}_h\times{\bf n_k}+\Phi({\bf e}_h)\times {\bf n}_k)\cdot\overline{(-\nabla\times{\bf u})} ~ds
\label{5.new2}
\end{split}
\end{equation}
and
\begin{equation}
\begin{split}
\int_\Omega \varepsilon{\bf e}_h \cdot \overline{\nabla p}~d\text{x}
& =\sum_{k=1}^N \int_{\Omega_k} \varepsilon_k{\bf e}_h \cdot \overline{\nabla p} ~d\text{x} \\
& =\sum_{k=1}^N ~\bigg(-\int_{\Omega_k}\varepsilon_k(div\,{\bf e}_h) \cdot \overline p ~dV +\int_{\partial\Omega_k}(\varepsilon_k{\bf e}_h\cdot{\bf n}_k)\cdot\overline{p} ~ds \bigg) \\
& =\sum_{k=1}^N ~\int_{\partial\Omega_k}(\varepsilon_k{\bf e}_h\cdot{\bf n}_k)\cdot \overline{p} ~ds
=\sum_{l<j} \int_{\Gamma_{lj}}\llbracket\varepsilon{\bf e}_h\cdot{\bf n}\rrbracket\cdot \overline{p}~ds.\label{5.new3}
\end{split}
\end{equation}

Define
$$ G = \sum_{k=1}^N \delta^{-1}||\nabla\times{\bf u}||_{0,\gamma_k}^2+\sum_{l<j} \bigg(\alpha^{-1}||\nabla\times{\bf u}||_{0,\Gamma_{lj}}^2 +\beta^{-1}\omega^2 ||{\bf u}||_{0,\Gamma_{lj}}^2  +\theta^{-1}||p||_{0,\Gamma_{lj}}^2 \bigg). $$

Substituting (\ref{5.new2})-(\ref{5.new3}) into (\ref{5.new4}), we get
\begin{equation*}
\begin{split}
|\varepsilon|_\text{\!min} ||{\bf e}_h||^2_{0,\!\Omega}\!
& \leq |\sum_{k=1}^N \int_{\gamma_k}(-{\bf e}_h\times{\bf n_k}+\Phi({\bf e}_h)\times {\bf n}_k)\cdot\overline{(-\nabla\times{\bf u})} ds \\
& \!+\!\sum_{l<j}\!\Big(\! \int_{\!\Gamma_{\!lj}}\!\llbracket{\bf e}_h\!\times\!{\bf n}\rrbracket\!\cdot\!\overline{(\nabla\!\times\!{\bf u})} ds \!+\!\int_{\!\Gamma_{\!lj}}\!\llbracket(\nabla\!\times\!{\bf e}_h)\times{\bf n}\rrbracket\!\cdot\!\overline{{\bf u}} ds \!+\!\int_{\!\Gamma_{\!lj}}\!\llbracket\varepsilon{\bf e}_h\!\cdot\!{\bf n}\rrbracket\!\cdot\!\overline{p} ds \!\Big)\!| \\
& \leq C \sum_{k=1}^N ||-{\bf e}_h\times{\bf n_k}+\Phi({\bf e}_h)\times {\bf n}_k||_{0,\gamma_k}\cdot ||\nabla\times{\bf u}||_{0,\gamma_k} \\
& +\sum_{l<j} \bigg( ||\llbracket{\bf e}_h\times{\bf n}\rrbracket||_{0,\Gamma_{lj}}\!\cdot\!||\nabla\times{\bf u}||_{0,\Gamma_{lj}}+||\llbracket\Psi({\bf e}_h)\times{\bf n}\rrbracket||_{0,\Gamma_{lj}}\!\cdot\! ||\omega{\bf u}||_{0,\Gamma_{lj}} \\
& +||\llbracket\varepsilon{\bf e}_h\cdot{\bf n}\rrbracket||_{0,\Gamma_{lj}} \cdot ||p||_{0,\Gamma_{lj}} \bigg).
\end{split}
\end{equation*}

Then, by Cauchy-Schwarz inequality, we further obtain
\begin{equation}
|\varepsilon|_{\text{min}} ||{\bf e}_h||^2_{0,\Omega} \leq ||{\bf e}_h||_{\bf V} \cdot G^\frac{1}{2}.\label{focus}
\end{equation}

In the following we estimate the $G$. Taking $\delta=\alpha$ in $G$ and using the trace inequality (see \cite[Thm.1.6.6]{trace1} and \cite[Thm.A.2]{trace2})
$$ ||u||_{0,\partial K}^2 \leq C(h_K^{-1}||u||_{0,K}^2+||u||_{0,K} |u|_{1,K}), \quad \forall ~u\in H^1(K) $$
and
$$ ||u||_{0,\partial K}^2 \leq C(h_K^{-1}||u||_{0,K}^2+h_K^{2\eta}|u|_{1/2+\eta,K}^2), \quad \forall ~u\in H^{1/2+\eta}(K) \quad(\eta > 0) $$
on each element $K\in{\mathcal T}_h$, leads to
\begin{equation*}
\begin{split}
G &\leq \sum_{k=1}^N \Big(\alpha^{-1}||\nabla\times{\bf u}||_{0,\partial\Omega_k}^2 +\beta^{-1}\omega^2 ||{\bf u}||_{0,\partial\Omega_k}^2  +\theta^{-1}||p||_{0,\partial\Omega_k}^2 \Big) \\
&\leq C \sum_{k=1}^N \Big( \alpha^{-1} h_k^{-1}||\nabla\times{\bf u}||_{0,\Omega_k}^2+\alpha^{-1}h_k^{2s}||\nabla\times{\bf u}||_{1/2+s,\Omega_k}^2+ \beta^{-1}\omega^2 h_k^{-1}||{\bf u}||_{0,\Omega_k}^2 \\
& \quad +\beta^{-1}\omega^2 h_k^{2s} ||{\bf u}||_{1/2+s,\Omega_k}^2 +\theta^{-1}h_k^{-1}||p||_{0,\Omega_k}^2+\theta^{-1}||p||_{0,\Omega_k} |p|_{1,\Omega_k} \Big) \\
& \leq C \big(\alpha^{-1} h^{-1}||\nabla\times{\bf u}||_{0,\Omega}^2+\alpha^{-1}h^{2s}||\nabla\times{\bf u}||_{1/2+s,\Omega}^2+ \beta^{-1} h^{-1}\omega^2||{\bf u}||_{0,\Omega}^2 \\
& \quad +\beta^{-1} h^{2s} \omega^2||{\bf u}||_{1/2+s,\Omega}^2 +\theta^{-1}h^{-1}||p||_{0,\Omega}^2+\theta^{-1}||p||_{0,\Omega} |p|_{1,\Omega} \big).
\end{split}
\end{equation*}

Then, by the stability estimates (\ref{esti})-(\ref{estiesti}) and Lemma \ref{wpes}, we can bound $G$ as follows
\begin{equation*}
\begin{split}
G
& \leq C \alpha^{-1}\big( h^{-1}+h^{2s}\omega^2 \big) \frac{|\varepsilon|^4_\text{max}}{|\varepsilon|^2_\text{min}} ||{\bf e}_h||_{0,\Omega}^2 \\
& \quad + C \beta^{-1}\big( h^{-1}+h^{2s}\omega^2 \big) \frac{|\varepsilon|^4_\text{max}}{|\varepsilon|^3_\text{min}} ||{\bf e}_h||_{0,\Omega}^2  + C \theta^{-1} h^{-1} \frac{|\varepsilon|^2_\text{max}}{|\varepsilon|^2_\text{min}} ||{\bf e}_h||_{0,\Omega}^2 \\
& \leq C |\varepsilon|^2_\text{min} \big( h^{-1}+h^{2s}\omega^2 \big) ||{\bf e}_h||_{0,\Omega}^2
\end{split}
\end{equation*}
with $\alpha= \frac{|\varepsilon|^4_\text{max}}{|\varepsilon|^4_\text{min}}, \quad \beta= \frac{|\varepsilon|^4_\text{max}}{|\varepsilon|^5_\text{min}}$ and $\theta= \frac{|\varepsilon|^2_\text{max}}{|\varepsilon|^4_\text{min}}$.

Thus we have
\begin{equation} \label{focus1}
G^{1/2}
\leq C |\varepsilon|_\text{min} \big( h^{-1/2}+h^{s}\omega \big) ||{\bf e}_h||_{0,\Omega}
\leq C |\varepsilon|_\text{min} h^{-1/2}\omega^{1/2} ||{\bf e}_h||_{0,\Omega}.
\end{equation}

Taking (\ref{focus1}) into (\ref{focus}), we obtain the estimate (\ref{most}).
$\hfill \Box$

For a domain $D\subset \Omega$, let $||\cdot||_{s,\omega,D}$ be the $\omega-$weighted Sobolev norm defined by
$$ ||v||_{s,\omega,D} = \sum_{j=0}^{s}\omega^{2(s-j)}|v|_{j,D}^2. $$

Let the mesh triangulation ${\mathcal T}_h$ satisfy the definition stated in \cite{pwdg} and set $\lambda=\text{min}_{K\in{\mathcal T}_h}\lambda_K$, where $\lambda_K$ is the positive parameter that depends only on the shape of an element $K$ of ${\mathcal T}_h$
. Let $r$ and $q$ be given positive integers satisfying $q\geq 2r+1$ and $q\geq 2(1+2^{1/\lambda})$. Let the number $p$ of plane-wave propagation directions be chosen as $p=(q+1)^2$.

The following approximation can be viewed as a version of Theorem 5.4 in \cite{pwdg}.

\begin{lemma} \label{define}
(\cite[Lemma 4.2]{long}) Let ${\bf E}$ be the solution of the equations (\ref{eq2 3}). Assume that ${\bf E}\in H^{r+1}(\text{curl,\,K})$. Then there exists $\bm\xi_{\bf E}\in {\bf V}_p({\mathcal T}_h)|_K$ such that
\begin{equation}
\begin{split}
||{\bf E}- \bm\xi_{\bf E}||_{j-1,\omega,K}
& \leq C \omega^{-2} (1 + (\omega h_K)^{q+j-r+8})e^{(\frac{7}{4}-\frac{3}{4}\rho)\omega h_K}h_K^{r+1-j} \\
& \quad\quad\quad \times
[q^{-\lambda_K(r+1-j)}+(\rho q)^{-\frac{q-3}{2}}pq] ||\nabla\times{\bf E}||_{r+1,\omega,K}
\label{4.new7}
\end{split}
\end{equation}
for every $1\leq j \leq r$.
\end{lemma}

Based on the previous lemmas, we can build a $L^2$ error estimate for the proposed method.

\begin{theorem} \label{imp}
Assume that the analytical solution ${\bf E}$ of the Maxwell problem (\ref{eq2 3}) belongs to $H^{r+1}(curl,\,\Omega)$ , $r\in \mathbb{N}$, $r \geq 3$. Let
the parameters in (\ref{minfun}) be defined by (\ref{parameters}) and let
${\bf E}_h\in {\bf V}_p({\mathcal T}_h)$ be the solution of (\ref{eq4 5}). Then, for large $p=(q+1)^2$, we have
\begin{equation}
||{\bf E}-{\bf E}_h||_{0,\Omega} \leq C \frac{|\varepsilon|^2_\text{max}}{|\varepsilon|^2_\text{min}} \omega^{-\frac{3}{2}} \frac{h^{r-2}}{q^{\lambda(r-\frac{3}{2})}} ||\nabla\times{\bf E}||_{r+1,\omega,\Omega},
\label{error-estimate1}
\end{equation}
where $C$ is a constant independent of $p$ but dependent on $\omega$ and $h$ only through the product $\omega h$ as an increasing function.
\end{theorem}

\paragraph*{Proof.}
Let ${\bf F}_h={\bm\xi}_{\bf E}\in {\bf V}_p({\mathcal T}_h)$ be defined by Lemma \ref{define}, then ${\bf E}_h-{\bf F}_h\in {\bf V}_p({\mathcal T}_h)$. Notice that ${\bf V}_p({\mathcal T}_h)\subset {\bf V}({\mathcal T}_h)$,
by the definition of ${\bf E}$ and ${\bf E}_h$, we have
$$ A({\bf E}-{\bf E}_h,{\bf E}_h-{\bf F}_h)=0. $$

Then we get by the Cauchy--Schwarz inequality
$$ ||{\bf E}-{\bf E}_h||^2_{\bf V}=A({\bf E}-{\bf E}_h,{\bf E}-{\bf F}_h)\leq||{\bf E}-{\bf E}_h||_{\bf V}\cdot||{\bf E}-{\bf F}_h||_{\bf V}, $$
which implies that
\begin{equation}
||{\bf E}-{\bf E}_h||_{\bf V}\leq||{\bf E}-{\bf F}_h||_{\bf V}.\label{new2}
\end{equation}

It suffices to estimate $||{\bf E}-{\bf F}_h||_{\bf V}$. For ease of notation, set ${\bm\epsilon}_h ={\bf E}-{\bf F}_h$. By the definition of the norm $||\cdot||_{\bf V}$, we get
\begin{equation}  \label{new3}
||{\bm\epsilon}_h||_{\bf V}^2 \leq C \frac{|\varepsilon|^4_\text{max}}{|\varepsilon|^4_\text{min}} \sum\limits_{k=1}^N\big(\int_{\partial\Omega_k}|\bm\epsilon_{h,k}|^2 ~ds
 +\omega^{-2}\int_{\partial\Omega_k}|(\nabla\times\bm\epsilon_{h,k})\times{\bf n}_k|^2 ~ds \big),
\end{equation}
where $\bm\epsilon_{h,k}=\bm\epsilon_h|_{\Omega_k}$. In an analogous way with the proof of Corollary 5.5 in \cite{pwdg}, we can prove that
\begin{equation*}
\int_{\partial\Omega_k}|\bm\epsilon_{h,k}|^2 ~ds
 + \omega^{-2}\int_{\partial\Omega_k}|(\nabla\times\bm\epsilon_{h,k})\times{\bf n}_k|^2 ~ds
 \leq C_1 \omega^{-4}\bigg(\frac{h}{q^{\lambda}}\bigg)^{2r-3}||\nabla\times{\bf E}||^2_{r+1,\omega,\Omega_k}.
\end{equation*}

Substituting the above inequality into (\ref{new3}) and combing (\ref{new2}), yields
\begin{equation} \label{new0}
||{\bf E}-{\bf E}_h||_{\bf V} \leq C \, \omega^{-2} \frac{|\varepsilon|^2_\text{max}}{|\varepsilon|^2_\text{min}} \bigg(\frac{h}{q^{\lambda}}\bigg)^{r-\frac{3}{2}}||\nabla\times{\bf E}||_{r+1,\omega,\Omega}.
\end{equation}

From the Lemma \ref{LL}, we already know
$$ ||{\bf E}-{\bf E}_h||_{0,\Omega}
\leq C \omega^{1/2}h^{-1/2} ||{\bf E}-{\bf E}_h||_{\bf V}. $$

Thus we get the desired estimate (\ref{error-estimate1}) by combining the above inequality with (\ref{new0}).

$\hfill \Box$

\section{The case with homogeneous material}

In this section we consider the special case that the coefficient $\varepsilon$ in (\ref{eq2 1}) is a constant on the whole domain $\Omega$. For this special case, which was considered in the most existing papers,
the proposed variational formula (\ref{eq3 6}) can be simplified.

According to the choices of the parameters (see (\ref{parameters})), we have $\alpha=\delta=1$, $\beta={1\over |\varepsilon|}$ and $\theta={1\over |\varepsilon|^2}$.
It is easy to see that the minimization functional (\ref{minfun}) becomes
\begin{equation*}
\begin{split}
J({\bf F})
& = \sum\limits_{k=1}^N\int _{\gamma_k}|-{\bf F}_k\times{\bf n}_k + \Phi({\bf F}_k)\times{\bf n}_k-{\bf g}|^2 ~ds \\
& +\sum\limits_{l<j} \Big(\int_{\Gamma_{lj}}|\llbracket{\bf F}\times{\bf n} \rrbracket|^2 ~ds +{1\over |\varepsilon|} \int_{\Gamma_{lj}} |\llbracket\Psi({\bf F})\times{\bf n} \rrbracket|^2 ~ds +\int_{\Gamma_{lj}} |\llbracket{\bf F}\cdot{\bf n} \rrbracket|^2 ~ds \Big) \\
& = \sum\limits_{k=1}^N\int _{\gamma_k}|-{\bf F}_k\times{\bf n}_k + \Phi({\bf F}_k)\times{\bf n}_k-{\bf g}|^2 ~ds \\
& +\sum\limits_{l<j} \Big(\int_{\Gamma_{lj}}|\llbracket{\bf F}\rrbracket|^2 ~ds +{1\over |\varepsilon|} \int_{\Gamma_{lj}} |\llbracket\Psi({\bf F})\times{\bf n} \rrbracket|^2 ~ds \Big),
\quad \forall {\bf F}\in {\bf V}({\mathcal T}_h).
\end{split}
\end{equation*}

Here we have used the equality $|\llbracket{\bf F}\times{\bf n} \rrbracket|^2+|\llbracket{\bf F}\cdot{\bf n} \rrbracket|^2 = |\llbracket{\bf F}\rrbracket|^2$ because of the
relation $\llbracket{\bf F}\rrbracket=\llbracket {\bf n}\times({\bf F}\times{\bf n})\rrbracket+\llbracket({\bf F}\cdot{\bf n}){\bf n}\rrbracket$ and the
orthogonality between $\llbracket{\bf n}\times({\bf F}\times{\bf n})\rrbracket$ and $\llbracket({\bf F}\cdot{\bf n}){\bf n}\rrbracket$.

Then the variational problem can be expressed as follows: Find ${\bf E}\in {\bf V}({\mathcal T}_h)$ such that
\begin{equation*}
\begin{split}
& \sum_{k=1}^N \int_{\gamma_k}(-{\bf E}_k\times{\bf n}_k +\Phi({\bf E}_k)\times{\bf n}_k-{\bf g})\cdot\overline{-{\bf F}_k\times{\bf n}_k+\Phi({\bf F}_k)\times{\bf n}_k} ~ds \\
& +\sum_{l<j}\big(\int_{\Gamma_{lj}}\llbracket{\bf E}\rrbracket\cdot\llbracket\overline{\bf F}\rrbracket ~ds +{1\over |\varepsilon|} \int_{\Gamma_{lj}} \llbracket \Psi({\bf E})\times{\bf n} \rrbracket\cdot\overline{ \llbracket \Psi({\bf F})\times{\bf n} \rrbracket} ~ds \big) =0,
\quad \forall ~{\bf F}\in {\bf V}({\mathcal T}_h).
\end{split}
\end{equation*}

Thus the sesquilinear $A(\cdot,\cdot)$ and the functional $(\bm\xi,~\cdot)$ are written as
\begin{equation*}
\begin{split}
A({\bf E},{\bf F})
& \!=\!\sum_{k=1}^N \int_{\gamma_k}(-{\bf E}_k\times{\bf n}_k +\Phi({\bf E}_k)\times{\bf n}_k)\cdot\overline{-{\bf F}_k\times{\bf n}_k +\Phi({\bf F}_k)\times{\bf n}_k} ~dS \\
& \!+\!\sum_{l<j}\! \big(\! \int_{\Gamma_{lj}}\!\llbracket{\bf E}\rrbracket\cdot\llbracket\overline{\bf F}\rrbracket ~ds +{1\over |\varepsilon|}\! \int_{\Gamma_{lj}}\!\llbracket\Psi({\bf E})\times{\bf n} \rrbracket\cdot\overline{ \llbracket\Psi({\bf F})\times{\bf n} \rrbracket} ~ds \big),
~ \forall ~{\bf F}\in {\bf V}({\mathcal T}_h).
\end{split}
\end{equation*}
and
\begin{equation*}
(\bm\xi,~{\bf F}) =
\sum_{k=1}^N\int _{\gamma_k}{\bf g}\cdot\overline{-{\bf F}_k\times{\bf n}_k+\Phi({\bf F}_k)\times{\bf n}_k} ~dS, \quad\quad \forall ~{\bf F}\in {\bf V}({\mathcal T}_h).
\end{equation*}

From the expression of $A(\cdot,\cdot)$, we know that the proposed method does not increase any cost of computation when $\varepsilon$ is a constant.
Moreover, for this case the error estimate in (\ref{error-estimate1}) becomes
\begin{equation}
||{\bf E}-{\bf E}_h||_{0,\Omega}
\leq C \omega^{-\frac{3}{2}} \frac{h^{r-2}}{q^{\lambda(r-\frac{3}{2})}} ||\nabla\times{\bf E}||_{r+1,\omega,\Omega}.\label{error-estimate2}
\end{equation}

\begin{remark} We emphasize that the $L^2$ error estimates (\ref{error-estimate1}) and (\ref{error-estimate2}) are almost the same as the $L^2$ error estimates of the plane wave approximations for three-dimensional Helmholtz 
equations in homogeneous media (there is an extra factor $\frac{|\varepsilon|^2_\text{max}}{|\varepsilon|^2_\text{min}}$ in (\ref{error-estimate1})).   
\end{remark}

\section{Numerical experiments}

For the examples discussed in this section, we adopt a uniform triangulation $\mathcal {T}_h$ for the domain $\Omega$ as follows: $\Omega$ is divided into small cubes of equal meshwidth, where $h$ is the length of the longest edge of the elements. As described in section 4.3, we choose the number of basis functions $p$ to be $p=(q+1)^2$ for all elements $\Omega_k$, where $q$ is a positive integer.

To measure the accuracy of the numerical solution, we introduce the following relative numerical error
$$ \text{err.}={||{\bf E}_{ex}-{\bf E}_h||_{L^2(\Omega)}\over{||{\bf E}_{ex}||_{L^2(\Omega)}}} \ $$
for the exact solution $ {\bf E}_{ex} \in L^2(\Omega)$, or
$$ \text{err.}=\sqrt{\sum_j^{Num}|{\bf E}_{ex,j}-{\bf E}_{h,j}|^2\over{\sum_j^{Num}|{\bf E}_{ex,j}|^2}}$$
for the exact solution ${\bf E}_{ex} \notin L^2(\Omega)$, where ${\bf E}_{ex,j}$ and ${\bf E}_{h,j}$ are the exact solution and numerical approximation, respectively, of the problem at vertices referred to by the subscript $j$.

For the examples tested in this section, we assume that $\mu=1$ and $\sigma=\sqrt{\mu/|\varepsilon|}$. We perform all computations on a Dell Precision T7610 workstation using MATLAB implementations.

\subsection{The case with homogeneous media}

\subsubsection{Electric dipole in free space for a smooth case}

In this part we consider the example tested in \cite{long} and \cite{hmm}.  In Maxwell system (\ref{eq2 3}) we choose \(\Omega=[-0.5,0.5]^3\) and set $\varepsilon=1+i$ ($\mu$ and $\sigma$ are defined above).
We compute the electric field due to an electric dipole source at the point ${\bf x}_0=(0.6,0.6,0.6)$. The dipole point source can be defined as the solution of a homogeneous Maxwell system (\ref{eq2 3}).
The exact solution of the problems is
\begin{equation}
{\bf E}_{\text{ex}}=-\text{i}\omega I \phi({\bf x},{\bf x}_0){\bf a} + \frac{I}{{\text{i}}\omega\varepsilon} \nabla(\nabla\phi\cdot{\bf a})
\end{equation}
with
$$ \phi({\bf x},{\bf x}_0) = \frac{\text{exp}(\text{i}\omega\sqrt{\varepsilon}|{\bf x}-{\bf x}_0|)}{4\pi|{\bf x}-{\bf x}_0|}.  $$

Then the boundary data is computed by
\begin{equation}
{\bf g}_{\text{ex}}=-{\bf E}_{ex}\times{\bf n}+\frac{\sigma}{i\omega\mu}((\nabla\times{\bf E}_{ex})\times{\bf n})\times{\bf n}.
\end{equation}

We would like to compare the relative \(L^2\) errors of the approximations generated by the proposed method and the PWLS method in \cite{long}. For convenience, throughout this section we use ``new PWLS" to represent
the method proposed in this paper and ``old PWLS" to represent the method introduced in \cite{long}. The comparison results are shown in Table \ref{tt1} and Table \ref{tt2}.


\begin{center}
\tabcaption{}
\label{tt1}
Errors of approximations with respect to $h$ for the case of $\omega=4\pi$.
\vskip 0.1in
\begin{tabular}{|c|c|c|c|c|} \hline
  & $h$ & $\frac{1}{4}$ & $\frac{1}{8}$ & $\frac{1}{16}$ \\ \hline
\multirow{2}*{ $p=16$ } & \text{old PWLS} & 0.5226 & 0.1469 & 0.0267 \\ \cline{2-5}
                        & \text{new PWLS}  & 0.5153 & 0.1349 & 0.0219 \\ \hline
\multirow{2}*{ $p=25$ } & \text{old PWLS} & 0.2686 & 0.0482 & 0.0061 \\ \cline{2-5}
                        & \text{new PWLS}  & 0.2630 & 0.0413 & 0.0035 \\ \hline
\end{tabular}
\end{center}

\begin{center}
\tabcaption{}
\label{tt2}
Errors of approximations with respect to $h$ for the case of $\omega=8\pi$.
\vskip 0.1in
\begin{tabular}{|c|c|c|c|c|} \hline
  & $h$ & $\frac{1}{4}$ & $\frac{1}{8}$ & $\frac{1}{16}$ \\ \hline
\multirow{2}*{ $p=16$ } & \text{old PWLS} & 0.8214 & 0.4206 & 0.0918 \\ \cline{2-5}
                        & \text{new PWLS}  & 0.8201 & 0.4153 & 0.0851 \\ \hline
\multirow{2}*{ $p=25$ } & \text{old PWLS} & 0.6372 & 0.1751 & 0.0214 \\ \cline{2-5}
                        & \text{new PWLS}  & 0.6359 & 0.1711 & 0.0180 \\ \hline
\end{tabular}
\end{center}

The results listed in Tables \ref{tt1} and \ref{tt2} indicate that the approximations generated by the new PWLS method have smaller errors than that generated by the old PWLS method in \cite{long}.

\subsubsection{Another example}
In this part, we consider another example in the homogeneous media. We also choose 
\(\Omega=[-0.5,0.5]^3\)
and assume that $\varepsilon=1+i$. The other quantities are kept the same. The analytical solution is given by
\begin{equation}
{\bf E}_{\text{ex}} = (\omega\sqrt{\varepsilon} xz ~\text{cos}(\omega\sqrt{\varepsilon}y),  -z~\text{sin}(\omega\sqrt{\varepsilon}y),   -~\text{cos}(\omega\sqrt{\varepsilon}x)^t.
\end{equation}

In Table \ref{tt4} and Table \ref{tt5} we give some comparisons for the \(L^2\) relative errors as in the last subsection.


\begin{center}
\tabcaption{}
\label{tt4}
Errors of approximations with respect to $h$ for the case of $\omega=4\pi$.
\vskip 0.1in
\begin{tabular}{|c|c|c|c|c|} \hline
  & $h$ & $\frac{1}{4}$ & $\frac{1}{8}$ & $\frac{1}{16}$ \\ \hline
\multirow{2}*{ $p=16$ } & \text{old PWLS} & 0.1573 & 0.0272 & 0.0038 \\ \cline{2-5}
                        & \text{new PWLS}  & 0.1532 & 0.0246 & 0.0029 \\ \hline
\multirow{2}*{ $p=25$ } & \text{old PWLS} & 0.0587 & 0.0057 & 0.0006 \\ \cline{2-5}
                        & \text{new PWLS}  & 0.0551 & 0.0044 & 0.0003 \\ \hline
\end{tabular}
\end{center}

\begin{center}
\tabcaption{}
\label{tt5}
Errors of approximations with respect to $h$ for the case of $\omega=8\pi$.
\vskip 0.1in
\begin{tabular}{|c|c|c|c|c|} \hline
  & $h$ & $\frac{1}{4}$ & $\frac{1}{8}$ & $\frac{1}{16}$ \\ \hline
\multirow{2}*{ $p=16$ } & \text{old PWLS} & 0.4792 & 0.1316 & 0.0218 \\ \cline{2-5}
                        & \text{new PWLS}  & 0.4732 & 0.1279 & 0.0201 \\ \hline
\multirow{2}*{ $p=25$ } & \text{old PWLS} & 0.3481 & 0.0531 & 0.0050 \\ \cline{2-5}
                        & \text{new PWLS}  & 0.3465 & 0.0513 & 0.0039 \\ \hline
\end{tabular}
\end{center}

Similarly to the first example, the results listed in Tables \ref{tt4} and \ref{tt5} indicate that the approximations generated by the new PWLS method are slightly more accurate.

\subsection{The case in layered media}

In this Subsection, we test two examples in layered media. For the case in layered media (i.e., $\varepsilon$ is not a constant), it is difficult to give an analytic solution of the homogeneous Maxwell system (\ref{eq2 3}).
In order to compute accuracies of the plane wave approximations generated by the proposed method, as usual we replace the analytic solution by a ¡°good¡± approximation generated by the standard finite element method with very fine grids.

\subsubsection{Electric dipole in free space for a smooth case}

Let the boundary data ${\bf g}$ in the equations (\ref{eq2 3}) be chosen as the same function ${\bf g}_{\text{ex}}$ given in Subsubsection 7.1.1. But, we choose $\varepsilon=1+i$ for the upper domain $z>0$ and $\varepsilon=2+2i$ for the region $z<0$ (the other quantities are kept the same as in Subsubsection 7.1.1). Since the parameter $\varepsilon$ is different from the one chosen in Subsubsection 7.1.1, the vector function ${\bf E}_{ex}$ in Subsubsection 7.1.1 is not the analytic solution of the current equations yet.

The $L^2$ errors of the resulting approximate solutions are shown in Table \ref{tt1p} and Table \ref{tt2p}.

\begin{center}
\tabcaption{}
\label{tt1p}
Errors of approximations with respect to $h$ for the case of $\omega=4\pi$.
\vskip 0.1in
\begin{tabular}{|c|c|c|c|c|} \hline
  & $h$ & $\frac{1}{4}$ & $\frac{1}{8}$ & $\frac{1}{16}$ \\ \hline
\multirow{2}*{ $p=16$ } & \text{old PWLS} & 0.5235 & 0.1475 & 0.0287 \\ \cline{2-5}
                        & \text{new PWLS}  & 0.5196 & 0.1399 & 0.0250 \\ \hline
\multirow{2}*{ $p=25$ } & \text{old PWLS} & 0.2697 & 0.0496 & 0.0128 \\ \cline{2-5}
                        & \text{new PWLS}  & 0.2662 & 0.0445 & 0.0115 \\ \hline
\end{tabular}
\end{center}

\begin{center}
\tabcaption{}
\label{tt2p}
Errors of approximations with respect to $h$ for the case of $\omega=8\pi$.
\vskip 0.1in
\begin{tabular}{|c|c|c|c|c|} \hline
  & $h$ & $\frac{1}{4}$ & $\frac{1}{8}$ & $\frac{1}{16}$ \\ \hline
\multirow{2}*{ $p=16$ } & \text{old PWLS} & 0.8220 & 0.4211 & 0.0921 \\ \cline{2-5}
                        & \text{new PWLS}  & 0.8207 & 0.4179 & 0.0878 \\ \hline
\multirow{2}*{ $p=25$ } & \text{old PWLS} & 0.6376 & 0.1767 & 0.0308 \\ \cline{2-5}
                        & \text{new PWLS}  & 0.6365 & 0.1745 & 0.0290 \\ \hline
\end{tabular}
\end{center}

The results listed in Tables \ref{tt1p} and \ref{tt2p} indicate that the approximations generated by our new PWLS method also have sightly smaller errors than that generated by the old PWLS method in \cite{long} for a layered medium.

\subsubsection{Another example}
In this part, we consider another example with more complicated structure of $\varepsilon$. We choose \(\Omega=[0,\,1]^3\) and define $\varepsilon$ as
\begin{equation*}
\varepsilon =
\begin{cases}
\frac{3}{2}+\frac{3}{2}i &  y<0.5 , z<0.5, \\
1+i &  y>0.5 , z<0.5, \\
1+i &  y<0.5 , z>0.5, \\
\frac{1}{2}+\frac{1}{2}i  &  y>0.5 , z>0.5.
\end{cases}
\end{equation*}
The other quantities are kept the same.. The analytical solution is given by
\begin{equation}
{\bf g}_{\text{ex}}=(x, \,y, \,x+y+z)^t.
\end{equation}

In the following tables we show the $L^2$ errors of the resulting approximate solutions.

\begin{center}
\tabcaption{}
\label{tt3p}
Errors of approximations with respect to $h$ for the case of $\omega=4\pi$.
\vskip 0.1in
\begin{tabular}{|c|c|c|c|c|} \hline
  & $h$ & $\frac{1}{4}$ & $\frac{1}{8}$ & $\frac{1}{16}$ \\ \hline
\multirow{2}*{ $p=16$ } & \text{old PWLS} & 0.1517 & 0.0452 & 0.0139 \\ \cline{2-5}
                        & \text{new PWLS}  & 0.1514 & 0.0435 & 0.0129 \\ \hline
\multirow{2}*{ $p=25$ } & \text{old PWLS} & 0.0724 & 0.0200 & 0.0061 \\ \cline{2-5}
                        & \text{new PWLS}  & 0.0704 & 0.0189 & 0.0058 \\ \hline
\end{tabular}
\end{center}

\begin{center}
\tabcaption{}
\label{tt4p}
Errors of approximations with respect to $h$ for the case of $\omega=8\pi$.
\vskip 0.1in
\begin{tabular}{|c|c|c|c|c|} \hline
  & $h$ & $\frac{1}{4}$ & $\frac{1}{8}$ & $\frac{1}{16}$ \\ \hline
\multirow{2}*{ $p=16$ } & \text{old PWLS} & 0.4148 & 0.1328 & 0.0375 \\ \cline{2-5}
                        & \text{new PWLS}  & 0.4132 & 0.1318 & 0.0369 \\ \hline
\multirow{2}*{ $p=25$ } & \text{old PWLS} & 0.2572 & 0.0554 & 0.0138 \\ \cline{2-5}
                        & \text{new PWLS}  & 0.2569 & 0.0540 & 0.0125 \\ \hline
\end{tabular}
\end{center}

The results listed in Tables \ref{tt3p} and \ref{tt4p} indicate that the advantage of the new PWLS method is kept even for a complicated layered medium.

\bibliographystyle{siamplain}

\end{document}